\documentclass[reqno]{amsart}
\usepackage{hyperref}\usepackage{latexsym}

\def\open#1{\setbox0=\hbox{$#1$}
\baselineskip = 0pt
\vbox{\hbox{\hspace*{0.4 \wd0}\tiny $\circ$}\hbox{$#1$}}
\baselineskip = 11pt\!}

\newtheorem{theorem}{Theorem}[section]
\newtheorem{lemma}[theorem]{Lemma}

\newtheorem{example}[theorem]{Example}
\newtheorem{definition}[theorem]{Definition}
\newcommand{\bt}{\begin{theorem}}\newcommand{\et}{\end{theorem}}
\newcommand{\blem}{\begin{lemma}}\newcommand{\elem}{\end{lemma}}
\newcommand{\bex}{\begin{example}}\newcommand{\eex}{\end{example}}
\newcommand{\bd}{\begin{definition}}\newcommand{\ed}{\end{definition}}

\newcommand{\eps}{\varepsilon}

\newcommand{\R}{\mathbb R}
\newcommand{\N}{\mathbb N}

\newcommand{\gs}{\ensuremath{{\mathcal G}^s} }

\newcommand{\ns}{\ensuremath{{\mathcal N}} }

\newcommand{\comp}{\subset\subset}

\newcommand{\beast}{\begin{eqnarray*}}
\newcommand{\eeast}{\end{eqnarray*}}

\newcommand{\al}{\alpha}

\newcommand{\Om}{\Omega}

\newcommand{\vphi}{\varphi}

\newcommand{\D}{{\mathcal D}}

\newcommand{\pa}{\partial}
 \newcommand{\G}{{\mathcal G}}
 
\newcommand{\CC}{{\mathcal C}}

\newcommand{\nn}{\nonumber}
\newcommand{\beq}{ \begin{equation} }\newcommand{\eeq}{\end{equation} }
\newcommand{\bea}{\begin{eqnarray}}\newcommand{\eea}{\end{eqnarray}}
\newcommand{\beas}{\begin{eqnarray*}}\newcommand{\eeas}{\end{eqnarray*}}
\newcommand{\beqs}{\begin{equation*}}\newcommand{\eeqs}{\end{equation*}}

\newcommand{\ben}{\begin{enumerate}}\newcommand{\een}{\end{enumerate}}
\newcommand{\ba}{\begin{array}}\newcommand{\ea}{\end{array}}
\newcommand{\brem}{\begin{thr} {\bf Remark. }\rm}
\newcommand{\ethi}{\end{thr}}

\newcommand{\id}[1]{\stackrel{\circ}{#1}\hspace*{-3pt}}

\newcommand{\supp}{\mathop{\mathrm{supp}}}

\def\fn{\open{f}\,}

\newcommand{\Gg}{{\mathcal G}_g}
\newcommand{\Eg}{{\mathcal E}^g_M}
\newcommand{\Ng}{{\mathcal N}_g}
\newcommand{\Ggg}{{\mathcal G}_{\tilde g}}
\newcommand{\Egg}{{\mathcal E}^{\tilde g}_M}
\newcommand{\Ngg}{{\mathcal N}_{\tilde g}}

\begin{document}

\title{Generalized solutions of the Vlasov-Poisson system with singular data}

\author[I. Kmit]{Irina Kmit}
\address{Institute for Applied Problems of Mechanics and Mathematics,
Ukrainian Academy of Sciences, Naukova St. 3b, 79060 Lviv, Ukraine}
\email{kmit@ov.litech.net}

\author[M. Kunzinger]{Michael Kunzinger}
\address{Fakult\"at f\"ur Mathematik, Universit\"at  Wien\\
Nordbergstrasse 15\\ 1090 Wien\\ Austria}
\email{Michael.Kunzinger@univie.ac.at}
\urladdr{http://www.mat.univie.ac.at/\~{}mike/}

\author[R. Steinbauer]{Roland Steinbauer}
\address{Fakult\"at f\"ur Mathematik, Universit\"at Wien\\
Nordbergstrasse 15\\ 1090 Wien\\ Austria}
\email{Roland.Steinbauer@univie.ac.at}
\urladdr{http://www.mat.univie.ac.at/\~{}stein/}

\thanks{Supported by the Austrian Science Fund's projects Z-50, Y-237
and P16742, and an OEAD-grant of I. Kmit.}

\begin{abstract}
We study spherically symmetric solutions of the Vlasov-Poisson system in the
context of algebras of generalized functions. This allows to model highly
concentrated initial configurations and provides a consistent setting for studying
singular limits of the system. The proof of unique solvability in 
our approach depends on new stability properties of
the system with respect to perturbations.
\end{abstract}
\maketitle
\section{Introduction}\label{sec:intro}
In kinetic theory one often considers collisionless ensembles of
classical particles which interact only by fields which they create
collectively. This situation is commonly referred to as the mean field
limit of a many particle system. More precisely such ensembles are described by
a phase-space distribution function $f:\R\times\R^3\times\R^3\to\R_0^+$, where
$\int_D f(t,x,v)dx\,dv$ gives the number of particles which at time $t$ have
their position $x$ and velocity $v$ in the region $D$ of phase-space $\R^6$.
The Vlasov equation expresses
the fact that $f$ is constant along particle paths---which is a direct
consequence of the absence of collisions---and reads
\begin{equation*}\pa_t f(t,x,v)+\pa_xf(t,x,v)\,v+\pa_vf(t,x,v)\,F(t,x)=0,\end{equation*} 
where $F$ is some force, which will emerge via some field equation
with its source given by the spatial particle density $\rho(t,x):=\int_{\R^3}f(t,x,v)\,dv$.
In the case of non-relativistic gravitational
or electrostatic fields the corresponding system of partial differential
equations is the Vlasov-Poisson system, in the case of relativistic electrodynamics
it is the Vlasov-Maxwell system and in the case of general relativistic gravity
the Vlasov-Einstein system. Such systems have been studied extensively in the literature; for an
overview see~\cite{glassey,pbull}.

The system which is best understood is the Vlasov-Poisson
system where $\rho$ acts as a source term for the Poisson equation.
For this system global-in-time classical solutions for general (compactly supported) initial data have been established in
\cite{Pf,LP,Sch}; for a review see \cite{Rrev,Rrev2}. The existence-theory for the other systems
mentioned above is not equally well understood. For the relativistic Vlasov-Maxwell system
global classical solutions are so far only known in special cases (e.g.\ \cite{GlStr,GlSch2,R1}),
while global weak solutions for general data were obtained in \cite{DL}. An investigation 
of the Vlasov-Einstein system was initiated in \cite{RR1}.

It is also interesting to note that the Vlasov-Poisson system possesses as (formal) singular limit
cases the Euler-Poisson system with pressure zero and the classical n-body problem, for both of which
in general no global-in-time solutions exist.
More precisely, in the first case if one considers a phase-space density function which 
is concentrated in $v$-space, i.e., $f(t,x,v)=\rho(t,x)\delta(v-w(t,x))$, 
where $\delta$ denotes the Dirac $\delta$-function and $w$ is a velocity field,  
then $f$ formally solves the Vlasov-Poisson system iff $(f,w)$ solves the pressure-less Euler-Poisson system. 
Similarly a density function $f$ concentrated in position and velocity space, i.e., $f(t,x,v)=
\sum_{k=1}^N\delta(x-x_k(t))\delta(v-v_k(t))$, formally solves the Vlasov-Poisson system
iff $x_k,\ v_k$ solve the n-body problem. One main problem here is---of course---the use of distributions
in the context of nonlinear equations. Only few rigorous results relating (approximating sequences of)
such concentrated solutions of the Vlasov-Poisson system to the solutions of the respective limit systems
have been achieved; see \cite{sandor} for the first case and \cite{neunzert} for the second.
The main interest of course would be to use the far better existence-theory in the case of the
Vlasov-Poisson system to learn something about the solutions of the related systems, e.g., in the
context of shell crossing singularities in the case of the pressure-less Euler-Poisson system.

In this work we propose the use of algebras of generalized functions (in the sense of J.F. Colombeau
\cite{col1,col2})
to study these singular limits of the Vlasov-Poisson system. As a first step we prove an existence
and uniqueness result for singular solutions, i.e., solutions concentrated either in position-space 
or in momentum-space (or both), to the spherically symmetric Vlasov-Poisson system in a suitable 
algebra of generalized functions, where the latter provides us with a consistent
framework for treating singular i.e., distributional solutions of nonlinear PDEs.
The fundamental strategy of solving PDEs with singular initial data in the setting of algebras of generalized functions
is regularization of singularities by convolution with a mollifier depending on a regularization
parameter $\eps$ and first solving the equation for fixed $\eps$ using existence theory in
the smooth setting. Proving existence and uniqueness of generalized solutions then
amounts to deriving asymptotic estimates with respect to the regularization parameter. 
This process may alternatively be seen as uncovering new asymptotic stability results of smooth 
solutions to the system under perturbations of the initial data, which in our view is of independent
interest. For a general discussion of applications of Colombeau theory to PDEs see \cite{MObook}.
Recent investigations into linear PDEs in this framework can be found in \cite{HO, HOP}.

We organize our presentation in the following way. In section \ref{sec:prelim} we collect some
well-known facts on the spherically symmetric Vlasov-Poisson system which will be used later on
and recall the basic definitions of generalized function algebras in the sense of J.F. Colombeau.
Our main results are stated and proved in section \ref{sec:main}. Finally we collect some facts on
solutions of the Poisson equation in this setting of generalized functions in an appendix.
 
Although our notation is mostly standard or self-explaining
we explicitly mention the following conventions:
For a function $h=h(t,x,v)$ or $h=h(t,x)$ we denote for given $t$
by $h(t)$ the corresponding function of the remaining variables.
By $\|\,.\,\|_p$ we denote the usual $L^p$-norm for $p\in[1,\infty]$.
The index $c$ in function spaces refers to compactly supported
functions. Constants denoted by $C$ may change their value from line to line 
but never depend on $\eps$.

\section{Preliminaries}\label{sec:prelim}

We start by collecting some preliminaries (for a comprehensive presentation including
full proofs see \cite{Rrev2}) from the existence-theory of the (spherically symmetric) 
Vlasov-Poisson system, which from now on we shall abbreviate by (VP)
\begin{eqnarray}
 \pa_tf+v\pa_xf-\pa_xu\pa_vf&=&0\label{eq:v}\\
 \Delta u&=&4\pi\gamma\rho\label{eq:p}\\
 \rho&=&\int_{\R^3}f\,dv
\end{eqnarray}
where $\gamma=\pm 1$. We suppose the following initial resp.\ boundary conditions
\begin{eqnarray}
 f(0,x,v)&=&\fn(x,v)\geq 0 \in \mathcal{C}^\infty_c(\R^6)\label{eq:ic}\\
 \lim_{|x|\to\infty}u(t,x)&=&0.
\end{eqnarray}
We shall often combine position and velocity into a single variable $z=(x,v)$ and denote
by $Z(s)=Z(s,t,z)=(X(s,t,z),V(s,t,z))$ the solutions of the characteristic system of (\ref{eq:v}), 
\begin{eqnarray*}
        \dot X(s) &=& V(s)\\
        \dot V(s) &=& -\pa_x u(s,X(s)) 
\end{eqnarray*}
with initial condition $Z(t,t,z)=z$. The solution
of the Vlasov equation is then given by $f(t,z)=\fn(Z(0,t,z))$, hence all $L^p$-norms ($1\leq p
\leq\infty$) of $f$ are constant in time, as is the $L^1$-norm of $\rho$, i.e., the mass, which 
will be denoted by $M$. Clearly $f(t)$ is compactly supported and we denote its velocity support 
by $P(t):=\sup\{|v|:\ (x,v)\in\supp(f(t))\}$. 

We shall call a function $g:\R^6=\R^3\times\R^3\to \R$ spherically symmetric if for all $A\in\mathrm{SO}(3)$
\begin{equation}\label{eq:ssym} g(Ax,Av)=g(x,v). \end{equation}
It is well known that in case the initial value $\fn$ of (VP) is spherically 
symmetric the respective solution $f(t)$ will also have this property.
Moreover, the spatial density $\rho(t)$ will be spherically symmetric 
(in the usual sense on $\R^3$---we shall denote it hence by $\rho(t,r)$, where
$r=|x|$) and the Poisson equation simplifies to
\[
\Delta u(t,r)=\frac{1}{r^2}\Big(r^2 u'(t,r)\Big)'=4\pi\gamma\rho(t,r). 
\]
By a slight abuse of notation, in what follows we will use $u(t,x)$ and
$u(t,r)$ interchangeably.
In addition to the usual key-estimate on the solution of the
Poisson equation with compactly supported source term $\rho(t)$, i.e.,
\begin{equation}\label{eq:u'}
 ||\pa_x u(t)||_\infty\leq C||\rho(t)||_1^{1/3}\,||\rho(t)||_\infty^{2/3}
\end{equation}
in the present setup we also obtain the estimates
\begin{eqnarray}
 |\pa_xu(t,r)|&\leq&\frac{M}{r^2}     \label{eq:sym2}\quad \mbox{and}\\ 
 ||\pa_x^{\alpha+e_i+e_j} u(t)||_\infty
  &\leq& C||\pa_x^\alpha\rho(t)||_\infty\qquad(\alpha\in\N_0^3,\ i,j=1,2,3).\label{eq:sym1}
\end{eqnarray}
Combining equations (\ref{eq:u'}) and (\ref{eq:sym2}) one obtains for all $r>0$
\begin{equation}\label{10}
 |\pa_x u(t,r)|\leq C\min\Big(\frac{1}{r^2},P(t)^2\Big).
\end{equation}
Note that the latter estimate together with the fact that from
$\xi\in\CC^2([0,t])$, $g\in L^1(\R)$ and $|\ddot\xi(s)|\leq g(\xi(s))\ 
\forall 0\leq s\leq t$ it follows that $|\dot\xi(t)-\dot\xi(0)|\leq 2\sqrt{2}||g||_1^{1/2}$ yields 
boundedness of $P(t)$ (cf.\ \cite{Rrev2}, proof of Thm.\ 1.4). Hence global existence of 
solutions follows by the standard continuation criterion. 

We now turn to algebras of generalized functions in the sense of J.F. Colombeau \cite{col1,col2}. These are differential 
algebras containing the vector space of distributions $\D'$ as a subspace and $\CC^\infty$ as a subalgebra while 
displaying maximal consistency with respect to classical analysis according to L. Schwartz' impossibility result \cite{sch-imp}. 
The main ingredient of the construction is regularization of distributions by nets
of smooth functions and asymptotic estimates in terms of the regularization parameter $\eps\in(0,1]=:I$, which in our 
case will be $L^\infty$-estimates global in $z$ on compact time intervals. 
We shall work within the so-called special version of the theory and use \cite{book} as our main 
reference. Colombeau algebras are defined as quotients of the spaces of moderate modulo negligible nets
$(u_\eps)_\eps$ in some basic space ${\mathcal E}$. In the present case we use 
${\mathcal E}=\CC^\infty(\R^+_0\times\R^n)^I$ and the following estimates for moderateness and negligibility
(where $O$ denotes the Landau symbol).

\begin{eqnarray}
 \Egg(\R^+_0\times\R^n)&:=&\{(u_\eps)_\eps\in{\mathcal E}:\ 
                        \forall K\comp\R^+_0\ \forall\alpha\in\N_0^{n+1}\ \exists N\in\N:
  \nn\\\label{def:mod}
                  &&\hfill \sup_{(t,z)\in K\times\R^n}|\pa^\al u_\eps(t,z)|=O(\eps^{-N})\
                      (\mbox{as}\ \eps\to 0)\},\\
 \Ngg(\R^+_0\times\R^n)&:=&\{(u_\eps)_\eps\in{\mathcal E}:\ 
                  \forall K\comp\R^+_0\ \forall\alpha\in\N_0^{n+1}\ \forall m\in\N:\nn\\\label{def:neg}  
                  &&\hfill \sup_{(t,z)\in K\times\R^n}|\pa^\al u_\eps(t,z)|=O(\eps^{m})\
                      (\mbox{as}\ \eps\to 0)\},\\ 
 \Ggg(\R^+_0\times\R^n)&:=&\Egg(\R^+_0\times\R^n)/\Ngg(\R^+_0\times\R^n). \nn
\end{eqnarray}
The index $\tilde g$ in the above definitions signifies the global estimates with respect to $x$, $v$ (contrary to the
local estimates w.r.t.\ $t$).
Generalized functions shall be denoted by $u=[(u_\eps)_\eps]$, meaning that $u$ is the equivalence class
of the net $(u_\eps)_\eps$. In the following section we shall prove 
existence and uniqueness results for the spherically symmetric (VP)-system in this setting.
By a solution of a differential equation in a Colombeau algebra $\Ggg$ we mean an element $[(u_\eps)_\eps]$ of the algebra
such that each $u_\eps$ solves the equations up to an element of the ideal $\Ngg$. Roughly, establishing
solvability of a PDE in this setting therefore amounts to obtaining a classical solution $u_\eps$ for each $\eps$ and
proving moderateness of the resulting net $(u_\eps)_\eps$. Proving uniqueness amounts to showing that any two nets 
solving the equation up to an element of the ideal (with initial data differing by an element of
the respective ideal) necessarily belong to the same equivalence class in $\Ggg$ (cf., e.g., \cite{MObook}). 

We note the following fact which will be essential in the uniqueness-part of the proof of our main result
(and follows by an easy adaptation of the proof of \cite[Thm.\ 1.2.3]{book}): $u=[(u_\eps)_\eps]\in\Egg(\R^+_0\times\R^n)$ is negligible iff
$$
\forall K\comp\R^+_0\ \forall m\in\N:\ \sup_{(t,z)\in K\times\R^n}|u_\eps(t,z)|=O(\eps^{m})\,.
$$

We shall also need a suitable algebra of generalized functions 
containing the initial data $\fn$. To this end we consider nets 
$(u_\eps)_\eps\in\CC^\infty(\R^n)^I$ satisfying estimates of the form
\[
  \forall \alpha\in\N_0^n\ \exists N\in\N\ \mbox{(resp. } \forall m\in\N):
  \sup_{x\in\R^n}|\pa^\al u_\eps(x)|=O(\eps^{-N})\ \mbox{(resp. } O(\eps^m)).
\]
We denote the respective spaces by $\Eg,\ \Ng$ and $\Gg$ (cf.\ \cite{BO, case}).

A function $u$ in $\Gg(\R^n)$ is called compactly supported  if there exists a representative 
$(u_\eps)_\eps$ of $u$ and a compact set $L$ containing the 
supports of all $u_\eps$. In this case we call the representative $(u_\eps)_\eps$ compactly supported.
Note however, that since $\Gg$ is not a sheaf there is no well-defined notion of support for its elements
(see Example \ref{nosheafex} below).      

The space $\D'_{L^\infty}(\R^n)$ of bounded distributions (distributional
derivatives of bounded functions) can be embedded into $\Gg(\R^n)$ by the map
$$
w \mapsto [(w*\vphi_\eps)_\eps]
$$
where $\vphi$ is a rapidly decreasing function with unit integral and all 
higher order moments vanishing, and $\vphi_\eps(x)=\eps^{-n}\vphi(x/\eps)$. 
This embedding commutes with partial derivatives. Analogously, $\mathcal{C}^\infty(\R^+_0,\D'_{L^\infty}(\R^n))$
can be embedded into $\Ggg(\R^+_0\times\R^n)$ via convolution.

\section{Generalized solutions of the Vlasov-Poisson system}\label{sec:main}

In this section we will state and prove our main results, providing existence
and uniqueness of generalized solutions of the spherically symmetric Vlasov-Poisson 
system. We begin with a discussion of the relevant symmetry properties.

We will call a generalized function $g\in\Gg(\R^6=\R^3\times\R^3)$ spherically
symmetric if it possesses a representative $(g_\eps)_\eps$ that is spherically
symmetric in the sense of (\ref{eq:ssym}) for all $\eps$. Likewise we call a function
$g\in\Gg(\R^3)$ spherically symmetric if it possesses a representative $(g_\eps)_\eps$ 
that for fixed $\eps$ is spherically symmetric in the usual sense.

The following definition singles out classes of scales which can be used to measure
the `maximal degree of divergence' admissible in the initial data of (VP) to allow
for unique solvability in the Colombeau algebra:

\bd Let $p>0$. 
\begin{itemize}
 \item[(i)] By $\Sigma_p^{(1)}$ we denote the space of all scales $\sigma: I \to I$ satisfying 
 $\sigma(\eps)\to 0$ for $\eps\to 0$ and
 \[\sigma(\eps)^{-1} = O(|\log(\eps)|^{1/p})\qquad (\eps\to 0)\,.\]
 \item[(ii)]  By $\Sigma_p^{(2)}$ we denote the space of all scales $\sigma: I \to I$ satisfying 
$\sigma(\eps)\to 0$ for $\eps\to 0$ and
 \[
\forall C>0 \ 
\exp(C\sigma(\eps)^{-p}) = O(|\log(\eps)|)\qquad (\eps\to 0)\,.
\]

\end{itemize}
\ed

Then $\Sigma_p^{(i)}\subseteq \Sigma_q^{(i)}$ for $p\ge q$ and $i=1,2$.
Note that using any scaling $\sigma$ satisfying $\sigma(\eps)\to 0$ for $\eps\to 0$, a $\delta$-source can be  
viewed as the element $[(\vphi_{\sigma(\eps)})_\eps]$ of the Colombeau algebra.
Since obviously $\vphi_{\sigma_\eps}\to \delta$ in $\D'$ as $\eps\to 0$, any such delta net
is associated to the standard image $[(\vphi_\eps)_\eps]$ of the Dirac measure,
hence macroscopically indistinguishable from it (cf.\ \cite{col2, MObook, book}
for discussions of the concept of association and its effects on nonlinear modelling
of singularities).
After these preparations we may state our existence result.
\bt\label{thm} (Existence of generalized solutions)\\
 Let $\fn\in\Gg(\R^6)$ with a spherically symmetric, non-negative and compactly
 supported representative $(\fn_\eps)_\eps$ satisfying
 \begin{enumerate}
  \item\label{1} $||\fn_\eps||_1=M$ (the mass) for all $\eps$, and
  \item\label{2} there exists some $\sigma\in \Sigma_2^{(1)}$ such that 
  $\displaystyle ||\fn_\eps||_\infty\leq\frac{C}{\sigma(\eps)}$. 
  \end{enumerate}
 Then there exists a solution $(f,u)$ of (VP) in $\Ggg(\R^+_0\times\R^6)
 \times\Ggg(\R^+_0\times\R^3)$ with $f(0,x,v)=\fn(x,v)$ and
 $u$ vanishing at infinity (in the sense of Definition ~\ref{vai}).
\et

Uniqueness of generalized solutions needs stronger assumptions on the data.
We present two results; the first one requires a $\Sigma_2^{(2)}$-scale.

\bt\label{u1}(Uniqueness of generalized solutions in $\Ggg$)\\
 Let the assumptions of Theorem \ref{thm} be satisfied but strengthen 
 (ii) to
 \begin{itemize}
  \item[(ii')] there exists some $\sigma\in\Sigma_2^{(2)}$ such that 
   $\displaystyle ||\fn_\eps||_\infty\leq\frac{C}{\sigma(\eps)}$. 
 \end{itemize}
Then $(f,u)$ given in Theorem \ref{thm} is the unique solution of (VP)
with $f(0,x,v)=\fn(x,v)$, velocity support of $f$ bounded as in Lemma \ref{lem1} (ii)
and $u$ vanishing at infinity (in the sense of Definition ~\ref{vai}).
\et

If we adjust the algebra to the symmetry of our problem we can do without a 
$\Sigma_p^{(2)}$-scale. More precisely we change the basic space ${\mathcal E}$ in the 
definitions of $\Egg$, $\Ngg$ and $\Ggg$, respectively to 
${\mathcal E}^\circ:=\{(u_\eps)_\eps\in{\mathcal E}=\CC^\infty(\R^+_0\times\R^n)^I|\
\forall t\in\R^+_0\ \forall \eps\in I:\ u_\eps(t)\mbox{ is spherically symmetric }\}$,
where in case $n=6$ spherical symmetry is to be understood in the sense of (\ref{eq:ssym})
and in case $n=3$ in the usual sense. We denote the resulting algebra by 
${\mathcal G}^\circ_{\tilde g}$. Likewise in case of the algebra $\G_g$ we take nets
$(u_\eps)_\eps\in\CC^\infty(\R^n)^I$ such that $u_\eps$ is spherically symmetric, again
in the respective senses for $n=6$ and $n=3$. The resulting algebra is denoted by ${\mathcal G}^\circ_g$.
Now we may state. 

\bt\label{u2} (Uniqueness of generalized solutions in $\Ggg^\circ$)\\
 Let $\fn\in\G^\circ_g(\R^6)$ with a compactly
 supported and non-negative representative $(\fn_\eps)_\eps$ satisfying (i)
 in Theorem \ref{thm} and
 \begin{enumerate}
  \item[(ii'')] there exists some $\sigma\in \Sigma_{\frac{10}{3}}^{(1)}$ such that 
   $\displaystyle ||\pa_z^\alpha \fn_\eps||_\infty\leq\frac{C}{\sigma(\eps)^{1+|\alpha|}}$ 
   for $|\alpha|\leq 1$.   
 \end{enumerate}
Then there exists a unique solution $(f,u)$ of (VP) in $\Ggg^\circ(\R^+_0\times\R^6)
\times\Ggg^\circ(\R^+_0\times\R^3)$ with $f(0,x,v)=\fn(x,v)$,
velocity support of $f$ bounded as in Lemma \ref{lem1} (ii) and
$u$ vanishing at infinity. 
\et

Note that the assumptions in the above theorems in particular allow to model concentrated data 
which lead to the singular limits of the Vlasov-Poisson system described in the introduction.

To prepare the proof of Theorem \ref{thm} first note that for fixed $\eps$ the classical theory provides 
us with a unique solution $(f_\eps,u_\eps)$ in $\CC^\infty(\R^+_0,\R^6)\times\CC^\infty(\R^+_0\times\R^3)$ 
with initial data $f_\eps(0,z)=\fn_\eps(z)$ and $u_\eps(t)\to 0$ as $|x|\to\infty$. Moreover, the solution
will inherit the symmetry property of the data, that is $f_\eps(t)$, $\rho_\eps(t)$ as well as 
$u_\eps(t)$ will be spherically symmetric.

To prove the existence of generalized solutions we have to verify the moderateness estimates in 
(\ref{def:mod}). We split this task into two Lemmas collecting the necessary estimates.

\blem\label{lem1} (Zero order estimates)
 There exists $C>0$ such that for $\eps$ sufficiently small we have for all $t\in\R^+$ 
 \begin{enumerate}
  \item $\displaystyle ||f_\eps(t)||_\infty\leq\frac{C}{\sigma(\eps)}$, 
      $||f_\eps(t)||_1=||\rho_\eps(t)||_1=M$
  \item $\displaystyle P_\eps(t) \leq\frac{C}{\sigma(\eps)^{\frac{1}{3}}}$ 
  \item $\displaystyle ||u_\eps(t)||_\infty,\ ||\pa_x u_\eps(t)||_\infty
         \leq\frac{C}{\sigma(\eps)^{\frac{4}{3}}}$
  \item $\displaystyle ||\rho_\eps(t)||_\infty\leq\frac{C}{\sigma(\eps)^2}$.
 \end{enumerate}
Moreover, for any $T>0$ and $\eps$ sufficiently small
\begin{itemize}
\item[(v)] $\sup_{t\in [0,T]}||Z_\eps(t)||_\infty \leq\frac{C}{\sigma(\eps)^{\frac{1}{3}}}$.
\end{itemize}
\elem

\begin{proof}
(i) The $\mathrm{L}^\infty$-estimate follows easily since for $\eps$ fixed 
 by the smooth result (see Section 2) we have $f_\eps(t,z)=\fn_\eps(Z_\eps(0,s,z))$. 
 The $\mathrm{L}^1$-estimates are immediate from assumption (i) in the theorem.\\
(ii) We conclude from equation (\ref{eq:u'}) and (i)
 \beq\label{(*)}
   ||\pa_xu_\eps(t)||_\infty 
   \leq C||\rho_\eps(t)||_1^{\frac{1}{3}}||\rho_\eps(t)||_\infty^{\frac{2}{3}}
   \leq C\frac{P_\eps(t)^2}{\sigma(\eps)^{\frac{2}{3}}}.
 \eeq
 We set
 $$
 g_\eps(t,r):=
 \min\left\{\frac{1}{r^2},
     \Big(\frac{P_\eps(t)}{\sigma(\eps)^{\frac{1}{3}}}\Big)^2\right\}
 $$
 Note that $g_\eps(s,r)\leq g_\eps(t,r)$ for $s\leq t$ since $P_\eps$ is monotonically 
 increasing.
 Then combining $|\pa_xu_\eps(t,r)|\leq M/r^2$ with the above estimate we obtain
 from the characteristic equation
 \[
   |\ddot X^i_\eps(s)|=|\pa_{x_i}u_\eps(s,X_\eps(s))|\leq C g_\eps(t,|X^i_\eps(s)|)
 \]
for $s\leq t$ and $1\le i \le 3$. Therefore by the standard argument mentioned below equation (\ref{10}) we obtain
 \[
   |\dot X^i_\eps(t)-\dot X^i_\eps(0)|\leq2C\sqrt{2}||g_\eps(t)||_1^{\frac{1}{2}}
 \]
 and are left with calculating the $L^1$-norm of $g_\eps(t)$. We have
 \beas
   \int\limits_{-\infty}^\infty|g_\eps(t,r)|dr
   &=&2\int\limits_0^\infty|g_\eps(t,r)|dr
   \,\leq\,2C\int\limits_0^{\frac{\sigma(\eps)^{1/3}}{P_\eps(t)}}
                  \Big(\frac{P_\eps(t)}{\sigma(\eps)^{\frac{1}{3}}}\Big)^2dr
     +2C\int\limits_{\frac{\sigma(\eps)^{1/3}}{P_\eps(t)}}^\infty
       \frac{1}{r^2}dr\\
   &=&2C\frac{P_\eps(t)}{\sigma(\eps)^{\frac{1}{3}}}
      +2C\frac{P_\eps(t)}{\sigma(\eps)^{\frac{1}{3}}}
   \,\leq C\frac{P_\eps(t)}{\sigma(\eps)^{\frac{1}{3}}}.
 \eeas
 Thus we obtain 
 \beq\label{14}
  |\dot X^i_\eps(t)-\dot X^i_\eps(0)|\leq 
         C\frac{P_\eps(t)^\frac{1}{2}}{\sigma(\eps)^\frac{1}{6}} 
 \eeq
 and hence from the definition of $P_\eps$
 \[
   P_\eps(t)\leq\,\id P+CP_\eps(t)^{\frac{1}{2}}\sigma(\eps)^{-\frac{1}{6}},
 \]
 where $\id{P}$\ \hspace{0.5mm} bounds the diameter of the support of $\fn_\eps$. This in turn implies that
 $P_\eps(t)$ is bounded independent of $t$ for $\eps$ fixed, and that $P_\eps(t)\leq\frac{C}{\sigma(\eps)^{\frac{1}{3}}}$, 
 which together with (\ref{14}) gives (ii) and (v).

 We insert (ii) into (\ref{(*)}) to prove (iii),
 i.e.,
 \[
   ||\pa_xu_\eps(t)||_\infty
   \leq \frac{C}{\sigma(\eps)^{\frac{4}{3}}}.
 \]
 The estimate on $u_\eps(t)$ now follows easily by integration (taking into account that
 $|u_\eps(t,x)| = O(1/|x|)$), while for (iv) we note
 \[
  ||\rho_\eps(t)||_\infty\leq C||\fn||_\infty P_\eps(t)^3\leq\frac{C}{\sigma(\eps)^2}.
 \]
 \end{proof}

\blem \label{lem2}(Higher order $x,v$-estimates)
For all $\alpha\in\N_0^6$, all $\beta\in\N_0^3$ and all $T>0$ there exists $C>0$ such that
for $\eps$ sufficiently small and all $t\in [0,T]$ we have
 \begin{enumerate}
  \item $\displaystyle ||\pa_z^\alpha f_\eps(t)||_\infty           \leq  e^{C\sigma(\eps)^{-2}}$
  \item $\displaystyle ||\pa_x^\beta\rho_\eps(t)||_\infty          \leq  e^{C\sigma(\eps)^{-2}}$
  \item $\displaystyle ||\pa_z^\alpha Z_\eps(t)||_\infty           \leq  e^{C\sigma(\eps)^{-2}}$
  \item $\displaystyle ||\pa_x^{\beta+e_i+e_j}u_\eps(t)||_\infty   \leq  e^{C\sigma(\eps)^{-2}}$ $\forall i,j$.
 \end{enumerate}
\elem

Note that compared with the zeroth order estimates we have to use an exponential term in $\sigma$ to 
bound the respective expressions necessitating 
the use of the scale $\sigma$ in condition~\ref{2} in Theorem~\ref{thm}. However, this term, i.e., 
$\exp(\sigma(\eps)^{-2})$ suffices to bound 
derivatives of {\em any} order. In particular, higher order derivatives do not lead to higher order 
exponential terms which would cause our approach to fail.

\begin{proof}
  We prove the Lemma by induction on $|\alpha|$ and $|\beta|$.\\
  In the case $|\alpha|=|\beta|=0$ we have shown even stronger estimates on $f_\eps(t)$, $Z_\eps(t)$ and
  $\rho_\eps(t)$ already in Lemma~\ref{lem1}. The only remaining estimate is the one on $\pa_x^2u_\eps(t)$ 
  which follows from $||\pa_x^2u_\eps(t)||_\infty\leq C||\rho_\eps(t)||_\infty\leq C\sigma(\eps)^{-2}$. 
  
  To carry out the inductive step we assume the Lemma holds for $|\alpha|,\ |\beta|\leq n$. We have to infer the respective
  estimates for $|\alpha|=|\beta|=n+1$. We define 
  \[ \xi^{(\alpha)}_\eps(s):=\pa_z^\alpha X_\eps(s),\quad \eta^{(\alpha)}_\eps(s):=\pa_z^\alpha V_\eps(s).\]
  Using the characteristic system we obtain (for suitably chosen $i$)
  \begin{eqnarray*}
   \dot\xi^{(\alpha)}_\eps(s)&=&\frac{d}{ds}\,\pa_z^{\alpha} X_\eps(s)
                            \,=\,\pa_z^\alpha V_\eps(s)\,=\,\eta^{(\alpha)}_\eps(s)\\
   \dot\eta^{(\alpha)}_\eps(s)&=&\frac{d}{ds}\,\pa_z^\alpha V_\eps(s)
                               \,=\,-\pa_z^{\alpha}\Big(\pa_x u_\eps\big(s,X_\eps(s,t,z)\big)\Big)\\
                           &=&-\pa_z^{\alpha-e_i}\Big(\pa_x^2u_\eps\big(s,X_\eps(s,t,z)\big)\pa_{z_i}X_\eps(s,t,z)\Big) \\ 
                           &=&-\pa_x^2 u_\eps\big(s,X_\eps(s)\big)\xi^{(\alpha)}_\eps(s)\\
                           &&-\sum_{0<\gamma\leq\alpha-e_i}\binom{\alpha-e_i}{\gamma}\pa_z^\gamma
                              \big(\pa_x^2 u_\eps(s,X_\eps(s))\big)\pa_z^{\alpha-\gamma}X_\eps(s)                                
  \end{eqnarray*}  
  The last expression is a sum of products of terms of the form
   \[ \pa_x^\delta u_\eps(s,X_\eps(s))\ \mbox{with}\  |\delta|\leq n+2\] 
  for which we have $||\pa_x^\delta u_\eps(t)||_\infty\leq C||\pa_x^{\delta'}\rho_\eps(t)||_\infty
  \leq \exp(C\sigma(\eps)^{-2})$ by (\ref{eq:sym1}) and the induction hypothesis since $|\delta'|\leq n$, and
   \[(\pa_z^\nu X_\eps)^{\nu'}\,(\pa_z^\omega X_\eps)^{\omega'}\ \mbox{with}\ \max(|\nu|,|\omega|)\leq n\] 
  which by induction hypothesis is also bounded by $\exp(C\sigma(\eps)^{-2})$ on compact time intervals.
  
  So we find using (\ref{eq:sym1}) for $|\alpha|=0$ and Lemma \ref{lem1} (iv)
   \[ |\dot\eta_\eps^{(\alpha)}(s)|\leq|\pa_x^2 u_\eps\big(s,X_\eps(s)\big)|\,|\xi_\eps^{(\alpha)}(s)| + e^{C\sigma(\eps)^{-2}}
                                   \le C\sigma(\eps)^{-2}|\xi_\eps^{(\alpha)}(s)| + e^{C\sigma(\eps)^{-2}} . \]
  Hence summing up we obtain
   \[ |\dot\eta_\eps^{(\alpha)}(s)|+|\dot\xi_\eps^{(\alpha)}(s)|\leq
                             e^{C\sigma(\eps)^{-2}}
                             +C\sigma(\eps)^{-2}\big(|\eta_\eps^{(\alpha)}(s)|+|\xi_\eps^{(\alpha)}(s)|\big),\]
  which by Gronwall's Lemma gives
  \begin{equation}
   |\eta_\eps^{(\alpha)}(s)|,\, |\xi_\eps^{(\alpha)}(s)|\leq e^{C\sigma(\eps)^{-2}},
  \end{equation}         
  i.e., $|\pa_z^{\alpha}Z_\eps(s)|\leq \exp(C\sigma(\eps)^{-2})$ on $[0,T]$ for all $|\alpha|=n+1$, which is (iii).
 
  From here we obtain 
   \[||\pa_z^\alpha f_\eps(t)||_\infty\leq Ce^{C\sigma(\eps)^{-2}}\ \mbox{for all}\ |\alpha|=n+1\]
  since $\pa_z^\alpha(\fn_\eps(t,Z_\eps(0,t,z)))$ is a sum of products of certain 
  $\pa_z^\delta\fn_\eps(Z_\eps(0,t,z))$ with products of powers of derivatives of $Z_\eps(0,t,z)$, and we can use (iii) and the 
  moderateness of $\fn_\eps$. Thereby we have also shown (i).
  
  Item (ii) is now obvious using Lemma~\ref{lem1} (ii). 
  Finally, to prove (iv) we combine (ii) with (\ref{eq:sym1}).

\end{proof}

{\em Proof Theorem \ref{thm}.}\\
By the estimates of Lemma~\ref{lem1} and Lemma~\ref{lem2} above we obtain the necessary bounds on
$\sup_{(t,z)\in K\times\R^6}|\pa^\al_z f_\eps(t,z)|$ and 
$\sup_{(t,x)\in K\times\R^3}|\pa^\beta_x u_\eps(t,z)|$, where $K$ is a compact subset of $[0,\infty)$.
 
To obtain the estimates on $\pa_tf_\eps$ we plug the estimates established so far into the Vlasov 
equation (using the bounded velocity support of $f_\eps(t)$). From here the estimate on $\pa_t\rho_\eps$ and hence on $\pa_tu_\eps$ follows. Now differentiating the Vlasov equation we obtain the estimates on terms of the form $\pa_t\pa_z^\al f_\eps$. Higher order $\pa_t$- and mixed $(t,z)$-estimates of $f_\eps$ are obtained by successively differentiating Vlasov's equation and in turn imply the respective estimates on $u_\eps$. 

Moreover $u(t)=[(u_\eps(t))_\eps]$ is vanishing at infinity in the sense of Definition \ref{vai} since the support of $\Delta u_\eps(t)$ is bounded by $C+tP_\eps(t)$ and $P_\eps(t)\leq C\sigma(\eps)^{-1/3}$ by Lemma \ref{lem1} (ii).

This proves existence of solutions in $\Ggg(\R^+\times\R^6)\times\Ggg(\R^+\times\R^3)$ with $f(0,x,v)=\fn(x,v)$ and $u(t)$ vanishing at infinity. \hfill$\Box$
\medskip

{\em Proof Theorem \ref{u1}.}\\
We have to prove uniqueness of the solution obtained above. 
So assume $(f=[(f_\eps)_\eps],u=[(u_\eps)_\eps])$ is a solution as 
constructed above and let $(\tilde f=[(\tilde f_\eps)_\eps],\tilde u=[(\tilde u_\eps)_\eps])$ be another 
solution of the (VP) system with the same initial data (i.e., $\tilde f(0)=\fn$), 
$\tilde u$ vanishing at infinity in the sense of Definition \ref{vai} (with distinguished representative
$\tilde u_\eps$), and $\eps$-wise bounded (by $\tilde P_\eps(t)$, satisfying (ii) of Lemma \ref{lem1}) 
velocity support of $\tilde f_\eps(t)$. Proving uniqueness in our setting amounts to establishing 
that the differences $f_\eps - \tilde f_\eps$ and $u_\eps - \tilde u_\eps$ lie in the respective ideals.
We have
 \begin{eqnarray}\label{pvp}
  \pa_t\tilde f_\eps+v\pa_x\tilde f_\eps-\pa_x\tilde u_\eps\pa_v\tilde f_\eps&=&n_\eps\nn\\ 
  \Delta \tilde u_\eps&=&4\pi\gamma\int_{\R^3}\tilde f_\eps\, dv+n_\eps\\
   \tilde f_\eps(0)&=&\fn_\eps+n_\eps=:\tilde f_\eps^\circ,  \nn
 \end{eqnarray} 
 where $(n_\eps)_\eps$ denotes a ``generic'' (analogous to the ``generic'' constant $C$) element of the ideal
 which may denote different negligible quantities in each equation. 
 Denoting by $\tilde Z_\eps$ the characteristics of the above inhomogeneous Vlasov equation we obtain
 (cf.\ e.g., \cite{KRS}, appendix A):
 \bea\label{tildef} \tilde f_\eps(t,z)
        &=&\tilde f_\eps^\circ(\tilde Z_\eps(0,t,z))+\int_0^tn_\eps(s,\tilde Z_\eps(s,t,z))\,ds\nn\\
        &=&\tilde f_\eps^\circ(\tilde Z_\eps(0,t,z))+n_\eps(t,z)=\fn_\eps(\tilde Z_\eps(0,t,z))+n_\eps(t,z).
 \eea  
Consequently we may estimate the difference in the distribution functions using Lemma \ref{lem2} (i)
 \bea\label{fdiff}
   |f_\eps(t,z)-\tilde f_\eps(t,z)|
    &\leq&||\pa_z\fn_\eps||_\infty|Z_\eps(0,t,z)-\tilde Z_\eps(0,t,z)|+|n_\eps(t,z)|\nn\\
    &\leq& e^{C\sigma(\eps)^{-2}}|Z_\eps(0,t,z)-\tilde Z_\eps(0,t,z)|+|n_\eps(t,z)|.
 \eea

 For the characteristics we obtain (for $0\le s \le t$)
 \bea\label{tildechar}
  |X_\eps(s)-\tilde X_\eps(s)|&\leq&\int\limits_s^t|V_\eps(s')-\tilde V_\eps(s')|\, ds'\nn\\
  |V_\eps(s)-\tilde V_\eps(s)|&\leq&
         \int\limits_s^t|\pa_x u_\eps(s',X_\eps(s'))-\pa_x\tilde u_\eps(s',\tilde X_\eps(s'))|\, ds'\nn\\
         &\leq&\int\limits_s^t\Big(|\pa_x u_\eps(s',X_\eps(s'))-\pa_x u_\eps(s',\tilde X_\eps(s'))|\nn\\
         &&\qquad\qquad +|\pa_x u_\eps(s',\tilde X_\eps(s'))-\pa_x \tilde u_\eps(s',\tilde X_\eps(s'))|\Big)\,ds'\nn\\
         &\leq&\sup\limits_{s\leq s'\leq t}||\pa_x^2 u_\eps(s')||_\infty\,
               \int\limits_s^t|X_\eps(s')-\tilde X_\eps(s')|\,ds'\nn\\
         &&\qquad +\int\limits_s^t|\pa_x u_\eps(s',\tilde X_\eps(s'))-\pa_x \tilde u_\eps(s',\tilde X_\eps(s'))|\,ds.
 \eea
 Now we turn to the perturbed Poisson equation for $\tilde u_\eps$, i.e.,
 \beq\label{p} \Delta\tilde u_\eps(t,x)=4\pi\gamma\tilde\rho_\eps(t,x) + n_\eps(t,x).\eeq
 Since $\tilde u$ is strongly vanishing at infinity the right hand side in the above equation has its support
 in $B_{\eps^{-N}}(0)$ for some $N\geq 0$. Furthermore by assumption 
 $\tilde f_\eps(t)$ has its $v$-support contained in some
 $B_{\tilde P_\eps(t)}(0)$ and hence its $x$-support bounded by $\id R+\int_0^t\tilde P_\eps(s)ds$, where $\id R \ $ bounds 
 the $x$-support of $\tilde f_\eps(0)$. This implies that the support of $\tilde f_\eps(t)$ is bounded by some
 $B_{\tilde Q_\eps(t)}$ with $\tilde Q_\eps(t)\leq C(t)\tilde P_\eps(t)$ where $C(t)$ depends linearly on time.

 As a consequence $n_\eps(t)$ in equation (\ref{p}) above has its support also contained in
 some $B_{\eps^{-N}}(0)$. Therefore we may define $\tilde n_\eps(t,x):=\int n_\eps(t,y)/|x-y|\,dy$, which is clearly in the
 ideal and finally we have found a representative  $\bar u_\eps:=\tilde u_\eps- \tilde n_\eps$ of 
 $[(\tilde u_\eps)_\eps]$ that satisfies the non-perturbed Poisson equation with source $\tilde \rho_\eps$, i.e., 
 $\Delta \bar u_\eps=4\pi\gamma\tilde\rho_\eps$. This in turn implies $\Delta(u_\eps-\bar u_\eps)=4\pi\gamma
 (\rho_\eps-\tilde\rho_\eps)$ and using (\ref{eq:u'}) we write
 \beas
  &&||\pa_x(u_\eps-\bar u_\eps)(t)||_\infty\\
  && \qquad\leq
   C\ ||\int\limits_{\R^3}(f_\eps(t,.,v)-\tilde f_\eps(t,.,v))dv||_1^{1/3}\ 
   ||\int\limits_{\R^3}(f_\eps(t,.,v)-\tilde f_\eps(t,.,v))dv||_\infty^{2/3}.
 \eeas
 On estimating the $L^1$-norm above we use $\bar Q_\eps(t):=\max(Q_\eps(t),\tilde Q_\eps(t))$, where
 $Q_\eps(t)$ denotes the respective bound on the support of $f_\eps(t)$ and write
 \beas
  ||\int\limits_{\R^3}(f_\eps(t,.,v)-\tilde f_\eps(t,.,v))dy||_1
   &\leq&\int\limits_{B_{\bar Q_\eps(t)}}|f_\eps(t,z)-\tilde f_\eps(t,z)|dz\\
   &\leq&C\bar Q_\eps(t)^6||f_\eps(t)-\tilde f_\eps(t)||_\infty\\
   &\leq&\frac{C}{\sigma(\eps)^2}||f_\eps(t)-\tilde f_\eps(t)||_\infty,
 \eeas
 for $t\in [0,T]$, where we have used Lemma \ref{lem1} (ii) in the last step. So we find 
 \beq\label{udiff}
  ||\pa_x(u_\eps-\bar u_\eps)(t)||_\infty\leq \frac{C}{\sigma(\eps)^{4/3}}||f_\eps(t)-\tilde f_\eps(t)||_\infty.
 \eeq
 
 Now applying Gronwall's lemma to (\ref{tildechar}) we obtain for all $q$ 
 \bea\label{Zdiff}
  |Z_\eps(s)-\tilde Z_\eps(s)|
  &\leq& Ce^{\sup\limits_{s\leq s'\leq t}||\pa_x^2u_\eps(s')||_\infty}
         \int\limits_s^t|\pa_x u_\eps(s',\tilde X_\eps(s'))-\pa_x \tilde u_\eps(s',\tilde X_\eps(s'))|\,ds\nn\\
  &\leq& e^{C\sigma(\eps)^{-2}}\int\limits_s^t||\pa_xu_\eps(s')-\pa_x\bar u_\eps(s')||_\infty ds'+
         \int\limits_s^t||n_\eps(s')||_\infty ds'\nn \\
  &\leq& e^{C\sigma(\eps)^{-2}}\Big(\frac{1}{\sigma(\eps)^{4/3}}\int\limits_s^t||f_\eps(s')-\tilde f_\eps(s')||_\infty ds'
         +\eps^q\Big),
 \eea
 where we have again used Lemma \ref{lem1} and (\ref{udiff}) above.
 Finally we combine (\ref{fdiff}) with (\ref{Zdiff}) and use Gronwall's lemma for the second time to 
 obtain for all $q$
 \beas
  \sup_{t\in [0,T]}||f_\eps(t)-\tilde f_\eps(t)||_\infty\leq\exp\big(\sigma^{-4/3} e^{C\sigma(\eps)^{-2}}\big)\, \eps^q,
 \eeas
 and due to our assumptions on the scale we see that the difference of the distribution functions is in the 
 ideal. From here the respective estimates on the difference of the spatial densities and on $||u_\eps(t)-
 \tilde u_\eps(t)||_\infty$ follow easily.  $\hfill\Box$
\medskip

{\em Proof of Theorem \ref{u2}:}\\
As for {\em existence} just observe that (ii'') implies (ii) and that by classical theory
the solution inherits the respective symmetry properties of the data.

To prove {\em uniqueness} we assume that $(f=[(f_\eps)_\eps],u=[(u_\eps)_\eps])\in \G^\circ_{\tilde g}(
\R^+_0\times\R^6)\times\G^\circ_{\tilde g}(\R^+_0\times\R^3)$ is a solution as 
constructed above and $(\tilde f=[(\tilde f_\eps)_\eps],\tilde u=[(\tilde u_\eps)_\eps])$ is another
such solution with the same initial data, $\tilde u$ vanishing at infinity (with distinguished representative
$(\tilde u_\eps)_\eps$) and the velocity support 
of $\tilde f_\eps(t)$ bounded by $\tilde P_\eps(t)$ satisfying (ii) of Lemma \ref{lem1}. We follow the proof of Theorem \ref{u1} 
up to estimate (\ref{tildechar}) but now using spherical symmetry we provide a stronger
estimate on  
\bea\label{paudiff}
 \int\limits_s^t|\pa_x u_\eps(s',\tilde X_\eps(s'))-\pa_x \tilde u_\eps(s',\tilde X_\eps(s'))|\,ds'
  \hphantom{xxxxxxxxxxxxxxxxxxxxx}\nn\\
  \hphantom{xxxxxxxxxxxxxxxxx}
  \le 4\pi\int\limits_s^t\left|\frac{\tilde X_\eps(s')}{\tilde r_\eps^3(s')}\right|
   \int\limits_0^{\tilde r_\eps(s')}l^2\big|\rho_\eps(s',l)-\tilde\rho_\eps(s',l)\big|\, dl\,ds',
\eea
where $\tilde r_\eps$ denotes the modulus of $\tilde X_\eps$. Note that this formula does not hold unless we 
use the algebra $\Ggg^\circ$ since in general $\tilde\rho_\eps(t)$ will not be spherically symmetric
due to the non-symmetric perturbations in (\ref{pvp}). 
Estimating the difference of the spatial densities we find
using $\bar P_\eps(t)=\max(P_\eps(t),\tilde P_\eps(t))$ as well as (\ref{tildef})
\beas
  |\rho_\eps(t,r)-\tilde\rho_\eps(t,r)|
 &\le&\int\limits_{B_{\bar P_\eps(t)}}|f_\eps(t,r,v)-\tilde f_\eps(t,r,v)|\,dv\\
 &=&\int\limits_{B_{\bar P_\eps(t)}}|\fn_\eps\big(Z_\eps(0,t,z)\big)-\fn_\eps\big(\tilde Z_\eps(0,t,z)\big)|\,dv
       +|n_\eps(t,r)|\\
 &\leq& \frac{C}{\sigma(\eps)}\ ||\pa_z\fn_\eps||_\infty\,||Z_\eps(0,t,.)-\tilde Z_\eps(0,t,.)||_\infty+|n_\eps(t,r)|.
\eeas
Inserting this into (\ref{paudiff}) we obtain using Lemma~\ref{lem1} 
\beas
 &&\int\limits_s^t|\pa_x u_\eps(s',\tilde X_\eps(s'))-\pa_x \tilde u_\eps(s',\tilde X_\eps(s'))|\,ds\\
 &&\qquad\leq \frac{C}{\sigma(\eps)}\ ||\pa_z\fn_\eps||_\infty\sup\limits_{s\leq s'\leq t}|\tilde X_\eps(s')|
  \int\limits_s^t\big(||Z_\eps(0,s',.)-\tilde Z_\eps(0,s',.)||_\infty\\
 &&\hphantom{xxxxxxxxxxxxxxxxxxxxxxxxxxxxxxxxxxxx} +\sup\limits_{r\in\R} |n_\eps(s',r)|\big) ds'\\
 &&\qquad\leq C\sigma(\eps)^{-\frac{10}{3}}\int\limits_s^t||Z_\eps(0,s',.)-\tilde Z_\eps(0,s',.)||_\infty ds'+C\eps^q
\eeas
for $t\in [0,T]$ and all $q$. Now combining this with (\ref{tildechar}) we obtain
\beas
 ||Z_\eps(0,s,.)-\tilde Z_\eps(0,s,.)||_\infty
 &\leq& C\sigma(\eps)^{-\frac{10}{3}}\int\limits_s^t||Z_\eps(0,s',.)-\tilde Z_\eps(0,s',.)||_\infty ds'+C\eps^q,
\eeas
which by Gronwall's lemma gives
\[ \sup_{t\in [0,T]}||Z_\eps(0,t,.)-\tilde Z_\eps(0,t,.)||_\infty\leq\eps^q e^{C\sigma(\eps)^{-\frac{10}{3}}}.\]
Hence by our assumption on the scale the difference of the characteristics is negligible. 
By (\ref{tildef}) this immediately implies that $[(f_\eps)_\eps]=[(\tilde f_\eps)_\eps]$ and 
so the same holds true for the spatial density as well as for the potential. \hfill $\Box$

\begin{appendix}
    \section{Uniqueness for generalized solutions of the Poisson equation}

In this appendix we collect some facts on the Poisson equation 
within the framework of nonlinear generalized functions.  We focus on the question of uniqueness, 
presenting a solution concept providing the existence of unique generalized solutions 
subject to a boundary condition generalizing the classical condition $u\to 0$ $(|x|\to \infty)$. 
Throughout this appendix we assume that $n\ge 3$ and write the Poisson equation as $\Delta u = \rho$. 
Also, we denote the fundamental solution of the Laplace equation by $C_n /|x|^{n-2}$.

In addition to the algebra $\Gg(\R^n)$ used in our main results
we also treat the case of the standard (special) Colombeau algebra $\gs(\Om)$ (with $\Om\subseteq\R^n$) which is defined using estimates on compact subsets of $\Om$, i.e., 
\beas \forall\al\in\N_0^n\ \forall K\comp\Om\ \exists N\in\N\ (\mbox{resp.\ }\forall m\in\N):\\
\qquad\qquad \sup_{x\in K}|\pa^\al u_\eps(x)|=O(\eps^{-N})\ (\mbox{resp. }O(\eps^m)). 
\eeas

We begin with some preliminaries. 
Let $u\in\gs(\Om)$. Then $u$ has compact support (that is: $\exists K\comp\Om:$
$u|_{\Om\setminus K}=0$) if and only if there exists a representative
$(u_\eps)_\eps$ of $u$ and $L\comp\Om$ such that $\supp(u_\eps)\subseteq L$
for all $\eps>0$. In this case we say that $(u_\eps)_\eps$ itself has compact support.

Indeed for any compactly supported $u$ we may choose a cut off function $\chi\in\D(\Om)$ such that
$\chi\equiv 1$ on a neighborhood of the support of $u$. Then for any representative $(u_\eps)_\eps$ of $u$ we construct a new representative $(\chi u_\eps)_\eps$ which vanishes outside the support 
of $\chi$.

However, in general $L$ will properly contain the support of $u$ in its interior. Indeed
let $u=\iota(\delta)$ (with $\iota$ denoting the embedding of distributions into the
algebra of generalized functions) then there clearly exist representatives that vanish outside 
any compact neighborhood of the origin. On the other hand there is no 
representative which vanishes outside the support of $u$.

Next we note that any generalized function $u\in\gs(\R^n)$ has a 
representative $(u_\eps)_\eps$ which vanishes at infinity, i.e.,
$u_\eps(x)\to 0$ $(|x|\to\infty)\ \forall\eps$. Indeed take any representative
of $u$ and multiply it with an $\eps$-dependent cut-off function
$\chi_\eps$ which is equal to unity inside a ball of radius 
$1/(2\eps)$ and vanishes outside a ball of radius $1/\eps$.
Moreover we have the following warning example of non-uniqueness of 
generalized solutions to the Laplace equation.

\bex\label{ex:nonunique}
We consider $\Delta u=0$ in $\gs(\R^n)$. Clearly $u=0$ is a solution. 
On the other hand we construct a solution $\tilde u$ as follows:
Set $\tilde u_\eps=\chi_\eps$ with $\chi_\eps$ as above. Then $\tilde u_\eps$ 
vanishes at infinity, $[(\tilde u_\eps)]=1$ and $\Delta\tilde u=0$.

However, there does not exist a representative $(\hat u_\eps)_\eps$ of $\tilde u$ 
such that $\supp(\Delta\hat u_\eps)$ is contained in some ball of
radius $R$ for all $\eps$. Indeed suppose to the contrary that 
$\Delta\hat u_\eps=n_\eps\in\ns^s(\R^n)$ with $\supp(n_\eps)
\in B_R(0)$ for all $\eps$. Then by classical uniqueness we have that
$\hat u_\eps(x)=C_n\int\frac{n_\eps(y)}{|x-y|^{n-2}}dy$ and hence $(\hat u_\eps)_\eps$
is in the ideal which is not possible. 
\eex

This observation motivates the following definition securing uniqueness of solutions to the Poisson equation.

\bd\label{solc}
Let $\rho\in{\gs}(\R^n)$ be compactly supported. We call $u\in \gs(\R^n)$ a {\em solution of
the Poisson equation vanishing at infinity} if $\Delta u=\rho$ and
if there exists a representative $(u_\eps)_\eps$ of $u$ that satisfies
\begin{itemize}
\item[(i)] $\forall\eps>0:\, \lim\limits_{x\to\infty}u_\eps(x)=0$, and
\item[(ii)] $(\Delta u_\eps)_\eps$ is compactly supported.
\end{itemize}
\ed

We may now state the following result.

\bt\label{uni}
 Let $\rho\in{\gs}(\R^n)$ be compactly supported. Then there exists one and only one
 solution of the Poisson equation
 \[ \Delta u=\rho \] 
 vanishing at infinity.
\et

Note that the assumptions on $u$ in Definition
\ref{solc} are not redundant. Indeed the compact support of $\rho$ guarantees the
existence of a representative $(u_\eps)_\eps$ satisfying property (ii) and there
also exists a representative $(\tilde u_\eps)_\eps$ of $u$ which vanishes at infinity.
However, in general $u_\eps\not=\tilde u_\eps$ and uniqueness may fail as
is explicitly demonstrated by the example above.

\begin{proof}Existence: By the above we may choose a compactly supported
representative $(\rho^c_\eps)_\eps$ of $\rho$ and define 
$$
u_\eps(x):= C_n \int \frac{\rho^c_\eps(y)}{|x-y|^{n-2}}\,dy
$$ 
By the classical theory $u_\eps$ satisfies
both requirements stated in the theorem.

Uniqueness: Let $u$, $\tilde u$ be two solutions as above and choose representatives
$(u_\eps)_\eps$ and $(\tilde u_\eps)_\eps$ satisfying (i) and (ii) in Definition
\ref{solc}. From the second property we conclude that $\Delta(u_\eps-\tilde u_\eps)=n_\eps$ is
compactly supported. By the first property we have $u_\eps-\tilde u_\eps\to 0$
$(|x|\to\infty)$. Hence by the classical theory $(u_\eps-\tilde u_\eps)(x)= C_n\int
\frac{n_\eps(y)}{|x-y|^{n-2}}\,dy$ which obviously is in the ideal.
\end{proof}

We now turn to the ``global'' algebra $\Gg$. The basic difference between $\Gg$ and
$\G$ is that due to the global estimates defining it, $\Gg$ is not a sheaf:
\bex \label{nosheafex} Define $u_\eps^{(m)} \in \CC^\infty((-m,m))$ to be $1$ for $\eps>1/m$ and $\exp(-1/\eps)$
for $\eps\le 1/m$. Choose a partition of unity $(\chi_m)_{m\in \N}$ subordinate to $((-m,m))_{m\in \N}$
and set
$$
u_\eps(x) := \sum_{m=1}^\infty \chi_m(x)u_\eps^{(m)}(x)\,.
$$
Then $(u_\eps)_\eps \in \Eg(\R) \setminus \Ng(\R)$, so $u = [(u_\eps)_\eps]$ provides an example
of a nonzero element of $\Gg(\R)$ whose restriction to each $(-m,m)$ is zero. 
\eex

Clearly in this setting Example \ref{ex:nonunique}
does not work since here $(\tilde u_\eps)_\eps$ is not a representative of the function $1$. 
This opens the possibility of relaxing condition (ii) in Definition \ref{solc} which is necessary 
in the context of the (VP)-system since $\Delta u_\eps(t)$ as constructed in the proof of Theorem~\ref{thm} 
will not be compactly supported. On the other hand we have proved $P_\eps(t)\leq C\sigma(\eps)^{-1/3}$ 
in Lemma \ref{lem1} (ii). This motivates the following definition which will provide us with the solution 
concept used in our main results.
\bd\label{vai}
 Let $\rho\in{\Gg}(\R^n)$ be compactly supported. We call $u\in{\Gg}(\R^n)$ a {\em solution of
 the Poisson equation vanishing at infinity} if $\Delta u=\rho$ and
 if there exists a representative $(u_\eps)_\eps$ of $u$ that satisfies
 \begin{itemize}
  \item[(i)] $\forall\eps>0:\, \lim\limits_{x\to\infty}u_\eps(x)=0$, and
  \item[(ii)] $\supp(\Delta u_\eps)_\eps\subseteq B_{\eps^{-N}}(0)$ for some $N\geq 0$.
 \end{itemize}
\ed

Note that again conditions (i) and (ii) are not redundant. Indeed take $u_\eps$ with
$u_\eps=1$ on $B_{e^{1/\eps}}(0)$ and vanishing outside a ball of twice that radius. 
Then (i) clearly holds but $\Delta u_\eps\not=0$ near 
$|x|=e^{1/\eps}$. The desired result in this framework is

\bt Let $\rho\in{\Gg}(\R^n)$ be compactly supported. Then there exists one and only one
 solution of the Poisson equation
 \[ \Delta u=\rho \] 
 vanishing at infinity.
\et

\begin{proof}
 \emph{Existence} is proved as in Theorem~\ref{uni}. 

 To prove \emph{uniqueness} suppose we have two solutions $u$, $\tilde u$ in $\Gg(\R^n)$ vanishing at infinity. Let $(u_\eps)_\eps$ and $(\tilde u_\eps)_\eps$ be representatives according to Definition \ref{vai}. By condition (ii) we have $\Delta(u_\eps-\tilde u_\eps)=n_\eps$ with $(n_\eps)_\eps$ in the ideal and $\supp(n_\eps)\subseteq B_{\eps^{-N}}(0)$ for some $N$. So 
 \[ |u_\eps-\tilde u_\eps|(x)\le C_n\int_{B_{\eps^{-N}}(0)}\frac{|n_\eps(y)|}{|x-y|^{n-2}}\,dy
   \leq C\eps^m\int_0^{\eps^{-N}}r\,dr\qquad\forall m,
 \]  
 hence is in the ideal.
\end{proof}

\end{appendix}

{\em Acknowledgment:} We would like to thank Gerhard Rein for several helpful discussions on the subject.
Irina Kmit is thankful to the members of the DIANA group for their hospitality 
during her stay at Vienna University.


\begin{thebibliography}{99}

\bibitem{BO} Biagioni, H., Oberguggenberger, M.,
Generalized solutions to Burgers' equation. {\em J. Diff. Eqs.} {\bf 97} 263 - 287, (1992).

\bibitem{col1} Colombeau, J.F.,
New Generalized Functions and Multiplication of Distributions.
{\em (North Holland, Amsterdam, 1984).}

\bibitem{col2} Colombeau, J.F.,
Elementary Introduction to New Generalized Functions. 
{\em (North Holland, Amsterdam 1985).}

\bibitem{DL} DiPerna, R. J., Lions, P.-L.,
Global weak solutions of Vlasov-Maxwell systems.
{\em Comm. Pure Appl. Math.}  {\bf 42}, no.\ 6, 729--757 (1989).

\bibitem{glassey} Glassey, R.,
The Cauchy Problem in Kinetic Theory.
{\em(SIAM, Philadelphia, PA, 1996).}

\bibitem{GlSch2}
Glassey, R.,  Schaeffer, J.,
Global existence for the relativistic
Vlasov-Maxwell system with nearly neutral initial data.
{\em Commun.\ Math.\ Phys.}\ {\bf 119}, 353--384\ (1988).

\bibitem{GlStr}
Glassey, R.,  Strauss, W. A.,
Absence of shocks in an initially dilute collisionless plasma.
{\em Commun.\ Math.\ Phys.}\
{\bf 113}, 191--208 (1987).

\bibitem{book} Grosser, M., Kunzinger, M., Oberguggenberger, M., Steinbauer, R.,
Geometric Theory of Generalized Functions.
\emph{Mathematics and its Applications} {\bf 537} 
(Kluwer Academic Publishers, Dordrecht, 2001.)

\bibitem{HO} H\"ormann, G., Oberguggenberger, M.,
Elliptic regularity and solvability for partial differential equations with Colombeau coefficients,
{\em Electron. J. Differential Equations}, Vol. 2004, (14), 1--30 (2004).

\bibitem{HOP} H\"ormann, G., Oberguggenberger, M., Pilipovi\'c, S.,
Microlocal hypoellipticity of linear partial differential operators with generalized 
functions as coefficients,
{\em Trans. Amer. Math. Soc.}, {\bf 358}, 3363--3383 (2006).

\bibitem{KRS} Kunzinger, M., Rein, G., Steinbauer, R.,
On local solutions of the relativistic Vlasov-Klein-Gordon system.
{\em Electron. J. Differential Equations}, Vol. 2005, No. 01, 1-17, 2005.  

\bibitem{LP}
Lions, P.-L., Perthame, B.,
Propagation of moments and regularity for the 3-dimensional
Vlasov-Poisson system.
{\em Invent.\ Math.}\ {\bf 105}, 415--430 (1991).

\bibitem{neunzert} Neunzert, H.,
An introduction to the nonlinear {B}oltzmann-{V}lasov equation, in C. Cercignani (Ed.) 
{\em Kinetic theories and the Boltzmann equation},
Lecture Notes in Math.\ {\bf 1048}, 60--110, (Springer, Berlin, 1984).

\bibitem{case} Oberguggenberger, M.,
Case study of a nonlinear, nonconservative, non-strictly hyperbolic system. 
{\em Nonlinear Anal.} {\bf 19}, 53 - 79, (1992).

\bibitem{MObook} Oberguggenberger, M.,
Multiplication of distributions and applications to partial differential
equations, Pitman Research Notes in Mathematics {\bf 259}, Longman 1992.

\bibitem{pbull} Perthame, B.,
Mathematical tools for kinetic equations, {\em Bull.\ Amer.\ Math.\ Soc.}
{\bf 41}(2), 205 - 244, (2004).

\bibitem{Pf}
Pfaffelmoser, K.,
Global classical solutions of the Vlasov-Poisson system in three
dimensions for general initial data.
{\em J.\ Diff.\ Eqns.}\
{\bf 95}, 281--303 (1992).

\bibitem{Rrev}
Rein, G.,
Selfgravitating systems in Newtonian theory-the Vlasov-Poisson system , 
Mathematics of gravitation, Part 1 (Warsaw, 1996), 179--194, Banach Center Publ.,
{\bf 41}, Part I, Polish Acad.\ Sci.\, Warsaw, 1997.

\bibitem{Rrev2}
Rein, G.,
Collissionless kinetic euqations from astrophysics---The Vlasov-Poisson system.
{\em Handbook of Differential Equations},  Evolutionary Equations. Vol. {\bf3}. 
Eds. C.M. Dafermos and E. Feireisl, Elsevier (2007)

\bibitem{R1}
Rein, G.,
Generic global solutions of the relativistic Vlasov-Maxwell
system of plasma physics.
{\em Commun.\ Math.\ Phys.}\
{\bf 135}, 41--78 (1990).

\bibitem{RR1}
Rein, G., Rendall, A.~D.,
Global existence of solutions of the spherically symmetric
Vlasov-Einstein system with small initial data.
{\em Commun.\ Math.\ Phys.}\ {\bf 150}, 561--583 (1992).

\bibitem{sandor}Dietz C., Sandor, V.,
The hydrodynamical limit of the Vlasov-Poisson System,
Transport Theory Statist.\ Phys.\ {\bf 28}(5), 499--520, (1999).

\bibitem{Sch}
Schaeffer, J.,
Global existence of smooth solutions to the Vlasov-Poisson system
in three dimensions.
{\em Commun.\ Part.\ Diff.\ Eqns.}\
{\bf 16}, 1313--1335 (1991).

\bibitem{sch-imp} Schwartz, L.,
Sur l'impossibilit\'e de la multiplication des distributions.
{\em C. R. Acad. Sci. Paris} {\bf 239}, 847-848 (1954).

\end{thebibliography}
\end{document}